\documentclass[11pt]{article}
\usepackage{amsmath,amssymb}
\usepackage{amsthm}
\usepackage[margin=2.5cm]{geometry}
\usepackage{subfigure}
\usepackage{graphicx}
\usepackage{epsfig}
\usepackage{setspace}
\newtheorem{theorem}{Theorem}
\newtheorem{definition}{Definition}
\newtheorem{conjecture}{Conjecture}

\title{Permutation--invariant Niven numbers}
\author{Huiling Wu$^1$ and Senyue Lou$^{2,3}$\\
\small{\it $^1$Department of Mathematics, Lishui University, Lishui 323000, China}\\
\small{\it $^2$School of Physical Science and Technology, Ningbo University, Ningbo, 315211, China} \\
\small{\it $^3$Institute of Fundamental Physics and Quantum Technology, Ningbo University, Ningbo, 315211, China}}
\date{}

\begin{document}

\maketitle

\begin{abstract}
This paper introduces permutation-invariant Niven numbers--a novel class of Niven numbers where all digit permutations (with leading zeros automatically ignored) must retain the Niven property.
We demonstrate that there exist infinitely many such numbers and that their magnitude is unbounded. Furthermore, we present an exhaustive search method for identifying permutation--invariant Niven numbers.
\end{abstract}
\leftline{\bf MSC: \rm 11A07}
\leftline{\bf Keywords: \rm Niven numbers, permutation invariance, congruence equation}
\hskip.1in

\section{Introduction}

Niven numbers (Harshad numbers)--integers divisible by the sum of their digits--were first conceptualized by Niven in a 1977 lecture \cite{Niven}, where he explored their existence properties within base--10. The term ``Niven number (NN)" was later coined by Kennedy \cite{Kennedy82}, who formalized their study by establishing their natural density (zero) and constructing infinitely many such numbers. These seminal works catalyzed research into the combinatorial and arithmetic properties of digit--sum--dependent sequences, extending to generalizations across bases, gaps, additive representations, and asymptotic distributions. This paper surveys key developments stemming from these foundational contributions.

Since the pioneering works of Niven and Kennedy \cite{Niven,Kennedy82}, various researches have been made by many authors. Cooper and Kennedy prove the existence of arbitrarily long sequences of consecutive Niven numbers (NNs) \cite{Cooper}. The constructive method is elegant but computationally intensive for large \(k\). It laid the groundwork for subsequent density analyses. Grundman refines Cooper--Kennedy's theorem by establishing explicit bounds for the minimal starting point of \(k\)--consecutive Niven sequences. The proof leverages modular arithmetic efficiently, though the bounds are not asymptotically tight. The paper of De Koninck and Doyon provides a rigorous asymptotic formula \(N(x) \sim c \frac{x}{\log x}\) for the counting function of NNs \cite{Koninck}. The probabilistic model used is insightful, but the constant \(c\)'s dependence on base \(b\) warrants further exploration. Cai characterizes Niven numbers in bases 2 and 3 \cite{Cai}, revealing structural properties unique to low bases. The results are technically accessible but limited in generalizability to arbitrary bases. Fredricksen, Ionascu and Luca
solve the minimality problem for NNs with fixed digit sum \(s\) \cite{Fredricksen}. The combinatorial approach is novel, though computational constraints limit results to \(s \leq 500\).
Sanna presents a novel approach \cite{Sanna}, combining deep tools from analytic number theory with combinatorial digit sum manipulations and additive combinatorics to establish the additive basis property for Niven numbers. Wilson's paper provides a valuable constructive proof \cite{Wilson}, explicitly demonstrating the existence of exceptionally small sequences of 20 consecutive NNs.

Symmetry analysis is an effective method throughout the natural sciences, particularly in integrable systems due to the presence of infinitely many local \cite{Olver} and nonlocal \cite{Lou} symmetries.

In this paper, we study the invariance of the Niven property under digit permutations.

\section{Definitions and Preliminaries}

\begin{definition}
A permutation-invariant Niven number (PINN) is a positive integer for which every digit permutation (without leading zeros) yields a Niven number.
\end{definition}

For a $k$--digit number $A_k \equiv a_k a_{k-1} \cdots a_1$ (where $a_k$ is the most significant digit),
the set of all digit permutations forms
\begin{equation}
\mathcal{A}_k = \bigl\{ s(A_k) \bigm| s \in \mathcal{S}_k \bigr\},
\label{Ak}
\end{equation}
where $\mathcal{S}_k$ denotes the symmetric group on $k$ elements. This group is generated by transpositions
\begin{equation}
\mathcal{S}_k = \langle \tau_2, \tau_3, \dots, \tau_k \rangle, \quad
\tau_i = (1,i) \quad \text{for} \quad i = 2, \dots, k,
\label{tauk}
\end{equation}
where each generator $\tau_i$ acts on $A_k$ by swapping the digits at positions $1$ and $i$:
\begin{equation}
\tau_i (a_k \cdots a_i \cdots a_2 a_1) = a_k \cdots a_1 \cdots a_2 a_i.
\label{tau_action}
\end{equation}

\section{Research Questions}

Now the questions related to the permutation--invariant numbers are:

\begin{enumerate}
\item What is the smallest PINN?
\item Are all repdigits (e.g., 111, 222, ...) PINN?
\item Can a number containing the digit 0 be PINN?
\item Do infinitely many PINNs exist? If yes, what is their asymptotic density?
\item For fixed digit--length $k$, how many PINNs exist?
\item How can we algorithmically generate all PINNs?
\item What are necessary and sufficient conditions for such numbers?
\item What conditions must the digit sum $s$ satisfy to ensure that all permutations (without leading zeros) of a digit set are divisible by $s$?
\item How do these numbers relate to established OEIS sequences (e.g., NNs, repdigits)? Do subclasses (primes, palindromes) exist?
\item What open questions remain? Can this concept extend to other bases (e.g., binary)?
\end{enumerate}

\section{Basic Properties}

\subsection{Single-digit PINNs}

For the case of single--digit numbers, it is straightforward that all are PINNs. Specifically,
\begin{equation}
\mbox{\rm NN}_1 \equiv \{1, 2, 3, 4, 5, 6, 7, 8, 9\},
\label{NN1}
\end{equation}
where $\mbox{\rm NN}_1$ denotes the set of 1--digit PINNs. We note that 1 is the smallest element in this set.

\subsection{Repdigits as PINNs}

Regarding the second problem, we first observe that not all repdigits are Niven numbers (NNs), and consequently, not all repdigits are PINNs. However, if a repdigit is an NN, it trivially qualifies as a PINN since its digit permutations are identical. For example, any 3--digit repdigit $aaa$ with $a \in \mbox{\rm NN}_1$ is an NN because
\begin{equation}
aaa \equiv 0 \pmod{3a}.
\end{equation}

For notational convenience in subsequent analysis, we define $a_{(n)}$ as the $n$--digit repdigit formed by repeating digit $a$. Representative cases include:
\begin{align*}
1_{(3)} &= 111, \\
5_{(4)} &= 5555.
\end{align*}

It is established that infinitely many repdigit NNs (which are simultaneously PINNs) exist. These admit the general form
\begin{equation}
N = \frac{(10^k - 1)a}{9}, \quad a \in \mbox{\rm NN}_1, \quad 10^k \equiv 1 \pmod{9ka}.
\label{0}
\end{equation}

To explicitly demonstrate the existence of infinitely many such repdigit NNs/PINNs, we propose the following solution conjecture for the congruence condition in \eqref{0}:

\begin{conjecture}\label{conj1}
The repdigit NN equation \eqref{0} admits solutions parameterized as
\begin{eqnarray}
k&=&k_{\{n\alpha\beta\gamma\delta\}}=3^{n}\cdot m_0^{\alpha}\cdot m_1^{\beta}\cdot m_{21}^{\gamma_1}\cdot m_{22}^{\gamma_2}
\cdot m_{31}^{\delta_1}\cdot m_{32}^{\delta_2}\cdot m_{33}^{\delta_3}\cdot m_{34}^{\delta_4}\cdot m_{35}^{\delta_5}\nonumber\\
&&  \text{for} \quad \{\alpha\beta\gamma\delta\}\equiv \{ \alpha, \beta, \gamma_1, \gamma_2, \delta_1, \delta_2,\ \ldots,\ \delta_5\} \in \mathbb{N}_0 = \{0, 1, 2, \ldots\},
\label{knmpqr}
\end{eqnarray}
subject to the following minimal--index constraints,
\begin{align*}
& n \geq 0 \quad \text{when} \quad \{\alpha \beta \gamma \delta\}=\{0\}, \\
& n \geq 1 \quad \text{when} \quad \alpha \neq 0, \\
& n \geq 2 \quad \text{when} \quad \beta \neq 0, \\
& n \geq 3 \quad \text{when} \quad \gamma=\{\gamma_1,\ \gamma_2\} \neq \{0\}, \\
& n \geq 4 \quad \text{when} \quad \delta=\{\delta_1,\ \ldots,\ \delta_5\} \neq \{0\},
\end{align*}
where,
$m_0=37,\ m_1=333667,\ m_{21}=757,\ m_{22}=440334654777631,\ m_{31}=163,\ m_{32}=9937,\ m_{33}=2462401,\ m_{34}=676421558270641,$ and $ m_{35}=130654897808007778425046117.$
\end{conjecture}

The validity of Conjecture \ref{conj1} has been computationally verified up to our system's overflow threshold of $\sim 10^{10^8}$. Moreover, although the special case of Conjecture~\ref{conj1} for $\gamma_2 = \delta_4 = \delta_5 = 0$ can be rigorously proved, this case is not the topic of the present paper and will be treated elsewhere.

Under this conjecture, the repdigit numbers $a_{(k_{n\alpha\beta\gamma\delta})}$ with $a \in \mbox{\rm NN}_1$ are both NNs and PINNs since they satisfy:
\begin{align}
a_{(k_{n\alpha\beta\gamma\delta})} &= \frac{a}{9}(10^{k_{n\alpha\beta\gamma\delta}} - 1) \equiv 0 \pmod{a\cdot k_{n\alpha\beta\gamma\delta}}, \label{3na} \\
\text{subject to:} \quad & n \geq 0 \;\; \text{if} \;\; m=p=q=r=s=t=0, \nonumber \\
& n \geq 1 \;\; \text{if} \;\; \alpha \neq 0, \nonumber \\
& n \geq 2 \;\; \text{if} \;\; \beta \neq 0, \nonumber \\
& n \geq 3 \;\; \text{if} \;\; \{\gamma_1,\ \gamma_2\} \neq \{0\}, \nonumber \\
& n \geq 4 \;\; \text{if} \;\; \{\delta_1,\ldots,\delta_5\} \neq \{0\}, \nonumber
\end{align}
where $k_{n\alpha\beta\gamma\delta}$ is defined in \eqref{knmpqr}. The conclusion \eqref{3na} holds conditionally on the truth of Conjecture \ref{conj1}.

Though the expression \eqref{3na} with \eqref{knmpqr} has offered infinitely many repdigit NNs/PINNs, there still exist infinitely many other NNs/PINNs. We do not expand this topic further.

\subsection{Numbers Containing Zero}

For the third question, we adopt the mathematical convention where leading zeros in a number are omitted:
\begin{equation}
0_{(i-1)} a_i \ldots a_k = a_i \ldots a_k \quad \text{for} \quad a_i \neq 0,
\label{convention}
\end{equation}
with $0_{(i-1)}$ denoting $i-1$ leading zeros. Under this convention, PINNs may contain internal or trailing zeros.

As an example, consider numbers of the form $10_{(n_1)}20_{(n_2)}$ and $20_{(m_1)}10_{(m_2)}$ for nonnegative integers $n_1,\ n_2,\ m_1,\ m_2$. Since any permutation of their digits (after removing leading zeros) has digit sum $1 + 2 = 3$, and because
\begin{equation}
10_{(n_1)}20_{(n_2)}=0,\ 20_{(m_1)}10_{(m_2)}=0 \pmod{3} \label{nm12}
\end{equation}
(by the congruence property of digit sums modulo 3), all permutations yield Niven numbers. Thus, these numbers are PINNs.

\subsection{Infinitude and Density}

For the fourth question, the answer is affirmative because the arbitrariness of the integers $n,\ m,\ p,\ q,$ $ r,\ s,\ t,\ n_1,\ n_2,\ m_1$, and $m_2$ in \eqref{3na} and \eqref{nm12} demonstrates the existence of infinitely many PINNs. The remainder of this paper will present additional infinite families of PINNs.

The density of PINNs is defined as their asymptotic density (or natural density), given by the limit of the proportion of PINNs among the first $N$ positive integers as $N \to \infty$. Denoting $N_{\text{PINNs}}(N)$ and $N_{\text{NNs}}(N)$ as the counts of PINNs and NNs up to $N$ respectively, we have
\begin{equation}
d_{\text{PINNs}} = \lim_{N \to \infty} \frac{N_{\text{PINNs}}(N)}{N} \leq \lim_{N \to \infty} \frac{N_{\text{NNs}}(N)}{N} = d_{\text{NNs}} = 0, \label{d}
\end{equation}
implying that the density of PINNs is zero. The inequality in \eqref{d} follows because PINNs constitute a subset of NNs, while the final equality is a known result for NNs \cite{Kennedy82}.

\section{Exhaustive Search for Small PINNs}

Solving the fifth problem presents a moderate challenge. For conciseness, we explicitly list all PINNs only for the cases $k=2$ and $k=3$; more general PINNs will be discussed in subsequent problems.

\subsection{2-digit PINNs}

For $k=2$, verification confirms that the complete set of PINNs consists of
\begin{equation}
\mbox{\rm NN}_{2}\equiv \{10,\ 12,\ 18,\ 20,\ 21,\ 24,\ 27,\ 30,\ 36,\ 40,\ 42,\ 45,\ 48,\ 50,\ 54,\ 60,\ 63,\ 70,\ 72, \ 80,\ 81,\ 84,\ 90\}.\label{NN2}
\end{equation}

Using the convention \eqref{convention}, the union of $\mbox{\rm NN}_{1}$ and $\mbox{\rm NN}_{2}$ can be defined as
\begin{eqnarray}
&&\mbox{\rm PINN}_{2}\equiv \mbox{\rm NN}_{1} \cup \mbox{\rm NN}_{2} ={\cal{S}}_2{\cal{N}}_2,\label{PINN2}\\
&&{\cal{N}}_2=\{10,\ 12,\ 18,\ 20,\ 24,\ 27,\ 30,\ 36,\ 40,\ 45,\ 48,\ 50, 60,\ 70, \ 80,\ 90\}.\label{CalN}
\end{eqnarray}

The $\mathcal{S}_2$-invariance of the 2-digit PINNs defined in \eqref{CalN} follows directly from their construction, as expressed by the symmetry relation:
\begin{equation}
\mathcal{S}_2(\mathrm{PINN}_2) = \mathrm{PINN}_2. \label{Inv2}
\end{equation}

Furthermore, the inclusion hierarchy
\begin{equation}
\mathrm{PINN}_2 \supset \mathrm{PINN}_1 = \mathrm{NN}_1 \label{2c1}
\end{equation}
holds as a consequence of the notational convention established in \eqref{convention}.

\subsection{3-digit PINNs}

After detailed verification, we establish that the 3-digit PINNs can be expressed as:
\begin{align}
\mathrm{PINN}_3 &\equiv \mathcal{S}_3 \mathcal{N}_3
= \mathcal{S}_3 \big( \mathcal{N}_{31} \cup \mathcal{N}_{32} \cup \mathcal{N}_{33} \cup \mathcal{N}_{34} \big), \label{PINN3} \\
\mathcal{N}_{31} &= \{ 100,\, 200,\, 300,\, 400,\, 500,\, 600,\, 700,\, 800,\, 900 \}, \label{N31} \\
\mathcal{N}_{32} &= \{ 120,\, 180,\, 240,\, 270,\, 360,\, 450,\, 480 \}, \label{N32} \\
\mathcal{N}_{33} &= \{ 1_{(3)},\, 2_{(3)},\, 3_{(3)},\, 4_{(3)},\, 5_{(3)},\, 6_{(3)},\, 7_{(3)},\, 8_{(3)},\, 9_{(3)} \}, \label{N33} \\
\mathcal{N}_{34} &= \{ 117,\, 126,\, 135,\, 144,\, 225,\, 234,\, 288,\, 468 \}. \label{N34}
\end{align}

The $\mathcal{S}_3$--invariance of $\mathrm{PINN}_3$ is manifest:
\begin{equation}
\mathcal{S}_3 \, \mathrm{PINN}_3 = \mathrm{PINN}_3. \label{Inv3}
\end{equation}

Furthermore, when accounting for the convention \eqref{convention}, we immediately deduce the containment relation:
\begin{equation}
\mathrm{PINN}_3 \supset \mathrm{PINN}_2. \label{3c2}
\end{equation}

\section{Algorithmic Generation of PINNs}

We now address the sixth problem by introducing a general algorithm for generating all PINNs.

To find $k$--digit PINNs, we exploit their permutation invariance property through a two--stage search algorithm:

\textbf{Step 1.} Search for PINNs composed of $k$ nonzero digits.

\textbf{Step 2.} Consider augmenting $j$--digit PINNs ($j < k$) by inserting $k-j$ zeros such that the resulting numbers remain PINNs.

We now employ this two-stage algorithm to generate $k$--digit PINNs for $k = 4$ to $9$.

\subsection{4-digit PINNs}

For $k=4$, one can verify that the PINNs containing no zero digits consist of precisely the following 12 elements:
\begin{equation}
    \mathcal{N}_{45} = \{ 1116,\  1125,\  1134,\  1224,\  1233,\  2223,\  2448,\  2268,\  2466,\  3699,\  4446,\  6669 \}.
\label{N45}
\end{equation}

The next step demonstrates that appending one or more zeros to any element of $\mathrm{PINN}_{3}$ preserves the \text{Niven property}. Thus, the $\mathcal{S}_4$--invariant PINNs take the form:
\begin{equation}
    \mathrm{PINN}_{4} = \mathcal{S}_4 \mathcal{N}_4 = \mathcal{S}_3 \left( \mathcal{N}_{41} \cup \mathcal{N}_{42} \cup \mathcal{N}_{43} \cup \mathcal{N}_{44} \cup \mathcal{N}_{45} \right),
\label{PINN4}
\end{equation}
where $\mathcal{N}_{45}$ is given by \eqref{N45} and
\begin{align}
    \mathcal{N}_{41} &= \{ 10_{(3)},\  20_{(3)},\  30_{(3)},\  40_{(3)},\  50_{(3)},\  60_{(3)},\  70_{(3)},\  80_{(3)},\  900 \}, \label{N41} \\
    \mathcal{N}_{42} &= \{ 1200,\  1800,\  2400,\  2700,\  3600,\  4500,\  4800 \}, \label{N42} \\
    \mathcal{N}_{43} &= \{ 1_{(3)}0,\  2_{(3)}0,\  3_{(3)}0,\  4_{(3)}0,\  5_{(3)}0,\  6_{(3)}0,\  7_{(3)}0,\  8_{(3)}0,\  9_{(3)}0 \}, \label{N43} \\
    \mathcal{N}_{44} &= \{ 1170,\  1260,\  1350,\  1440,\  2250,\  2340,\  2880,\  4680 \}. \label{N44}
\end{align}

\subsection{Higher-digit PINNs}

Similarly, applying the two--step search algorithm to $k$--digit PINNs for $k = 5, 6, \ldots, 9$, we obtain:

\begin{align}
\mbox{\rm PINN}_{5}&\equiv {\cal{S}}_5{\cal{N}}_5={\cal{S}}_5({\cal{N}}_{51} \cup {\cal{N}}_{52} \cup {\cal{N}}_{53} \cup {\cal{N}}_{54} \cup {\cal{N}}_{55} \cup {\cal{N}}_{56}),\label{PINN5}\\
\mbox{\rm PINN}_{6}&\equiv {\cal{S}}_6{\cal{N}}_6={\cal{S}}_6({\cal{N}}_{61} \cup {\cal{N}}_{62} \cup {\cal{N}}_{63} \cup {\cal{N}}_{64} \cup {\cal{N}}_{65} \cup {\cal{N}}_{66} \cup {\cal{N}}_{67}),\label{PINN6}\\
\mbox{\rm PINN}_{7}&\equiv {\cal{S}}_7{\cal{N}}_7={\cal{S}}_7({\cal{N}}_{71} \cup {\cal{N}}_{72} \cup {\cal{N}}_{73} \cup {\cal{N}}_{74} \cup {\cal{N}}_{75} \cup {\cal{N}}_{76} \cup {\cal{N}}_{77} \cup {\cal{N}}_{78}),\label{PINN7}\\
\mbox{\rm PINN}_{8}&\equiv {\cal{S}}_8{\cal{N}}_8={\cal{S}}_8({\cal{N}}_{81} \cup {\cal{N}}_{82} \cup {\cal{N}}_{83} \cup {\cal{N}}_{84} \cup {\cal{N}}_{85} \cup {\cal{N}}_{86} \cup {\cal{N}}_{87} \cup {\cal{N}}_{88} \cup {\cal{N}}_{89}),\label{PINN8}\\
\mbox{\rm PINN}_{9}&\equiv {\cal{S}}_9{\cal{N}}_9={\cal{S}}_9({\cal{N}}_{91} \cup {\cal{N}}_{92} \cup {\cal{N}}_{93} \cup {\cal{N}}_{94} \cup {\cal{N}}_{95} \cup {\cal{N}}_{96} \cup {\cal{N}}_{97} \cup {\cal{N}}_{98} \cup {\cal{N}}_{99} \cup {\cal{N}}_{910}),\label{PINN9}
\end{align}
with the specific sets ${\cal{N}}_{ki}$ being defined as,
\begin{eqnarray*}
&&{\cal{N}}_{51}=\{10_{(4)},\ 20_{(4)},\ 30_{(4)},\ 40_{(4)},\ 50_{(4)}, 60_{(4)},\ 70_{(4)}, \ 80_{(4)},\ 90_{(4)}\},\nonumber\\
&&{\cal{N}}_{52}=\{120_{(3)},\ 180_{(3)},\ 240_{(3)},\ 270_{(3)},\ 360_{(3)},\ 450_{(3)},\ 480_{(3)}\},\nonumber\\
&&{\cal{N}}_{53}=\{1_{(3)}00,\ 2_{(3)}00,\ 3_{(3)}00,4_{(3)}00,\ 5_{(3)}00,\ 6_{(3)}00,\ 7_{(3)}00,\ 8_{(3)}00,\ 9_{(3)}00\},\nonumber\\
&&{\cal{N}}_{54}=\{11700,\ 12600,\ 13500,\ 14400,\ 22500,\ 23400,\ 28800,\ 46800\}, \nonumber\\
&&{\cal{N}}_{55}=\{1_{(3)}60, 11250, 11340, 12240, 12330, 2_{(3)}30, 24480, 22680, 24660, 36990, 4_{(3)}60, 6_{(3)}90\},\nonumber\\
&&{\cal{N}}_{56}=\{1_{(4)}5,1_{(3)}24,1_{(3)}33,11223,12_{(4)},2_{(3)}48,2_{(3)}66,22446,24_{(4)},3_{(3)}99,33669,36_{(4)},48_{(4)}
\},
\end{eqnarray*}
\begin{eqnarray*}
{\cal{N}}_{61}&=&\{10_{(5)},\ 20_{(5)},\ 30_{(5)},\ 40_{(5)},\ 50_{(5)}, 60_{(5)},\ 70_{(5)}, \ 80_{(5)},\ 90_{(5)}\},\nonumber\\
{\cal{N}}_{62}&=&\{120_{(4)},\ 180_{(4)},\ 240_{(4)},\ 270_{(4)},\ 360_{(4)},\ 450_{(4)},\ 480_{(4)}\},\nonumber\\
{\cal{N}}_{63}&=&\{1_{(3)}0_{(3)},\ 2_{(3)}0_{(3)},\ 3_{(3)}0_{(3)},4_{(3)}0_{(3)},\ 5_{(3)}0_{(3)},\ 6_{(3)}0_{(3)},\ 7_{(3)}0_{(3)},\ 8_{(3)}0_{(3)},\ 9_{(3)}0_{(3)}\},\nonumber\\
{\cal{N}}_{64}&=&\{1170_{(3)},\ 1260_{(3)},\ 1350_{(3)},\ 1440_{(3)},\ 2250_{(3)},\ 2340_{(3)},\ 2880_{(3)},\ 4680_{(3)}\}, \nonumber\\
{\cal{N}}_{65}&=&\{1_{(3)}600,\ 112500,\ 113400,\ 122400,\ 123300,\ 2_{(3)}300,\nonumber\\
 && 244800,\ 226800,\ 246600,\ 369900,\ 4_{(3)}600,\ 6_{(3)}900\},\nonumber\\
{\cal{N}}_{66}&=&\{1_{(4)}50,\ 1_{(3)}240,\ 1_{(3)}330,\ 112230,\ 12_{(4)}0,\ 2_{(3)}480,\nonumber\\
&&2_{(3)}660,\ 224460,\ 24_{(4)}0,\ 3_{(3)}990,\ 336690,\ 36_{(4)}0,\ 48_{(4)}0\},\nonumber\\
{\cal{N}}_{67}&=&\{1_{(5)}4,\ 1_{(4)}23,\ 1_{(3)}2_{(3)},\ 2_{(5)}8,\ 2_{(4)}46,\ 2_{(3)}4_{(3)},\ 3_{(4)}69,\ 3_{(3)}6_{(3)},\ 4_{(3)}8_{(3)}
\},
\end{eqnarray*}
\begin{eqnarray*}
{\cal{N}}_{71}&=&\{10_{(6)},\ 20_{(6)},\ 30_{(6)},\ 40_{(6)},\ 50_{(6)}, 60_{(6)},\ 70_{(6)}, \ 80_{(6)},\ 90_{(6)}\},\nonumber\\
{\cal{N}}_{72}&=&\{120_{(5)},\ 180_{(5)},\ 240_{(5)},\ 270_{(5)},\ 360_{(5)},\ 450_{(5)},\ 480_{(5)}\},\nonumber\\
{\cal{N}}_{73}&=&\{1_{(3)}0_{(4)},\ 2_{(3)}0_{(4)},\ 3_{(3)}0_{(4)},4_{(3)}0_{(4)},\ 5_{(3)}0_{(4)},\ 6_{(3)}0_{(4)},\ 7_{(3)}0_{(4)},\ 8_{(3)}0_{(4)},\ 9_{(3)}0_{(4)}\},\nonumber\\
{\cal{N}}_{74}&=&\{1170_{(4)},\ 1260_{(4)},\ 1350_{(4)},\ 1440_{(4)},\ 2250_{(4)},\ 2340_{(4)},\ 2880_{(4)},\ 4680_{(4)}\}, \nonumber\\
{\cal{N}}_{75}&=&\{1_{(3)}60_{(3)},\ 11250_{(3)},\ 11340_{(3)},\ 12240_{(3)},\ 12330_{(3)},\ 2_{(3)}30_{(3)},\nonumber\\
 && 24480_{(3)},\ 22680_{(3)},\ 24660_{(3)},\ 36990_{(3)},\ 4_{(3)}60_{(3)},\ 6_{(3)}90_{(3)}\},\nonumber\\
{\cal{N}}_{76}&=&\{1_{(4)}500,\ 1_{(3)}2400,\ 1_{(3)}3300,\ 1122300,\ 12_{(4)}00,\ 2_{(3)}4800,\nonumber\\
&&2_{(3)}6600,\ 2244600,\ 24_{(4)}00,\ 3_{(3)}9900,\ 3366900,\ 36_{(4)}00,\ 48_{(4)}00\},\nonumber\\
{\cal{N}}_{77}&=&\{1_{(5)}40,\ 1_{(4)}230,\ 1_{(3)}2_{(3)}0,\ 2_{(5)}80,\ 2_{(4)}460,\ 2_{(3)}4_{(3)}0,\ 3_{(4)}690,\ 3_{(3)}6_{(3)}0,\ 4_{(3)}8_{(3)}0
\},\nonumber\\
{\cal{N}}_{78}&=&\{1_{(6)}3,\ 1_{(5)}22,\ 2_{(6)}6,\ 2_{(5)}44,\ 3_{(6)}9,\ 3_{(5)}66,\ 4_{(5)}88
\},
\end{eqnarray*}
\begin{eqnarray*}
{\cal{N}}_{81}&=&\{10_{(7)},\ 20_{(7)},\ 30_{(7)},\ 40_{(7)},\ 50_{(7)}, 60_{(7)},\ 70_{(7)}, \ 80_{(7)},\ 90_{(7)}\},\nonumber\\
{\cal{N}}_{82}&=&\{120_{(6)},\ 180_{(6)},\ 240_{(6)},\ 270_{(6)},\ 360_{(6)},\ 450_{(6)},\ 480_{(6)}\},\nonumber\\
{\cal{N}}_{83}&=&\{1_{(3)}0_{(5)},\ 2_{(3)}0_{(5)},\ 3_{(3)}0_{(5)},4_{(3)}0_{(5)},\ 5_{(3)}0_{(5)},\ 6_{(3)}0_{(5)},\ 7_{(3)}0_{(5)},\ 8_{(3)}0_{(5)},\ 9_{(3)}0_{(5)}\},\nonumber\\
{\cal{N}}_{84}&=&\{1170_{(5)},\ 1260_{(5)},\ 1350_{(5)},\ 1440_{(5)},\ 2250_{(5)},\ 2340_{(5)},\ 2880_{(5)},\ 4680_{(5)}\}, \nonumber\\
{\cal{N}}_{85}&=&\{1_{(3)}60_{(4)},\ 11250_{(4)},\ 11340_{(4)},\ 12240_{(4)},\ 12330_{(4)},\ 2_{(3)}30_{(4)},\nonumber\\
 && 24480_{(4)},\ 22680_{(4)},\ 24660_{(4)},\ 36990_{(4)},\ 4_{(3)}60_{(4)},\ 6_{(3)}90_{(4)}\},\nonumber\\
{\cal{N}}_{86}&=&\{1_{(4)}50_{(3)},\ 1_{(3)}240_{(3)},\ 1_{(3)}330_{(3)},\ 112230_{(3)},\ 12_{(4)}0_{(3)},\ 2_{(3)}480_{(3)},\nonumber\\
&&2_{(3)}660_{(3)},\ 224460_{(3)},\ 24_{(4)}0_{(3)},\ 3_{(3)}990_{(3)},\ 336690_{(3)},\ 36_{(4)}0_{(3)},\ 48_{(4)}0_{(3)}\},\nonumber\\
{\cal{N}}_{87}&=&\{1_{(5)}400,\ 1_{(4)}2300,\ 1_{(3)}2_{(3)}00,\ 2_{(5)}800,\ 2_{(4)}4600,\nonumber\\ &&2_{(3)}4_{(3)}00,\ 3_{(4)}6900,\ 3_{(3)}6_{(3)}00,\ 4_{(3)}8_{(3)}00\},\nonumber\\
{\cal{N}}_{88}&=&\{1_{(6)}30,\ 1_{(5)}220,\ 2_{(6)}60,\ 2_{(5)}440,\ 3_{(6)}90,\ 3_{(5)}660,\ 4_{(5)}880
\},\nonumber\\
{\cal{N}}_{89}&=&\{1_{(7)}2,\ 2_{(7)}4,\ 3_{(7)}6,\ 4_{(7)}8\},
\end{eqnarray*}
\begin{eqnarray*}
{\cal{N}}_{91}&=&\{10_{(8)},\ 20_{(8)},\ 30_{(8)},\ 40_{(8)},\ 50_{(8)}, 60_{(8)},\ 70_{(8)}, \ 80_{(8)},\ 90_{(8)}\},\nonumber\\
{\cal{N}}_{92}&=&\{120_{(7)},\ 180_{(7)},\ 240_{(7)},\ 270_{(7)},\ 360_{(7)},\ 450_{(7)},\ 480_{(7)}\},\nonumber\\
{\cal{N}}_{93}&=&\{1_{(3)}0_{(6)},\ 2_{(3)}0_{(6)},\ 3_{(3)}0_{(6)},4_{(3)}0_{(6)},\ 5_{(3)}0_{(6)},\ 6_{(3)}0_{(6)},\ 7_{(3)}0_{(6)},\ 8_{(3)}0_{(6)},\ 9_{(3)}0_{(6)}\},\nonumber\\
{\cal{N}}_{94}&=&\{1170_{(6)},\ 1260_{(6)},\ 1350_{(6)},\ 1440_{(6)},\ 2250_{(6)},\ 2340_{(6)},\ 2880_{(6)},\ 4680_{(6)}\}, \nonumber\\
{\cal{N}}_{95}&=&\{1_{(3)}60_{(5)},\ 11250_{(5)},\ 11340_{(5)},\ 12240_{(5)},\ 12330_{(5)},\ 2_{(3)}30_{(5)},\nonumber\\
 && 24480_{(5)},\ 22680_{(5)},\ 24660_{(5)},\ 36990_{(5)},\ 4_{(3)}60_{(5)},\ 6_{(3)}90_{(5)}\},\nonumber\\
{\cal{N}}_{96}&=&\{1_{(4)}50_{(4)},\ 1_{(3)}240_{(4)},\ 1_{(3)}330_{(4)},\ 112230_{(4)},\ 12_{(4)}0_{(4)},\ 2_{(3)}480_{(4)},\nonumber\\
&&2_{(3)}660_{(4)},\ 224460_{(4)},\ 24_{(4)}0_{(4)},\ 3_{(3)}990_{(4)},\ 336690_{(4)},\ 36_{(4)}0_{(4)},\ 48_{(4)}0_{(4)}\},\nonumber\\
{\cal{N}}_{97}&=&\{1_{(5)}40_{(3)},\ 1_{(4)}230_{(3)},\ 1_{(3)}2_{(3)}0_{(3)},\ 2_{(5)}80_{(3)},\ 2_{(4)}460_{(3)},\nonumber\\ &&2_{(3)}4_{(3)}0_{(3)},\ 3_{(4)}690_{(3)},\ 3_{(3)}6_{(3)}0_{(3)},\ 4_{(3)}8_{(3)}0_{(3)}\},\nonumber\\
{\cal{N}}_{98}&=&\{1_{(6)}300,\ 1_{(5)}2200,\ 2_{(6)}600,\ 2_{(5)}4400,\ 3_{(6)}900,\ 3_{(5)}6600,\ 4_{(5)}8800
\},\nonumber\\
{\cal{N}}_{99}&=&\{1_{(7)}20,\ 2_{(7)}40,\ 3_{(7)}60,\ 4_{(7)}80\},\nonumber\\
{\cal{N}}_{910}&=&\{1_{(9)},\ 2_{(9)},\ 3_{(9)},\ 4_{(9)},\ 5_{(9)},\ 6_{(9)},\ 7_{(9)},\ 8_{(9)},\ 9_{(9)}\}.
\end{eqnarray*}

Apart from the repdigit PINNs such as those in \eqref{3na}, no $k$--digit PINNs without zero digits have been found for $k \geq 10$.

\section{Main Theorem}

\begin{theorem}\label{mainthm}
The $k$--digit integers with $k\geq10$
\begin{align}
\mbox{\rm PINN}_{k}&\equiv {\cal{S}}_k{\cal{N}}_k={\cal{S}}_k({\cal{N}}_{ka} \cup {\cal{N}}_{kb} \cup {\cal{N}}_{kc} \cup {\cal{N}}_{kd} \cup {\cal{N}}_{ke} \cup {\cal{N}}_{kf} \cup {\cal{N}}_{kg} \cup {\cal{N}}_{kh} \cup {\cal{N}}_{ki} \cup {\cal{N}}_{kj}),\label{PINNk}
\end{align}
are all PINNs, where ${\cal{N}}_{ka},\ {\cal{N}}_{kb},\ {\cal{N}}_{kc},\ {\cal{N}}_{kd},\ {\cal{N}}_{ke},\ {\cal{N}}_{kf},\ {\cal{N}}_{kg},\ {\cal{N}}_{kh},\ {\cal{N}}_{ki},$ and $ {\cal{N}}_{kj}$ are defined by
\begin{eqnarray*}
{\cal{N}}_{ka}&=&\{10_{k_1},\ 20_{k_1},\ 30_{k_1},\ 40_{k_1},\ 50_{k_1}, 60_{k_1},\ 70_{k_1}, \ 80_{k_1},\ 90_{k_1}\},
\\
{\cal{N}}_{kb}&=&\{120_{k_2},\ 180_{k_2},\ 240_{k_2},\ 270_{k_2},\ 360_{k_2},\ 450_{k_2},\ 480_{k_2}\},
\\
{\cal{N}}_{kc}&=&\{1_{(3)}0_{k_3},\ 2_{(3)}0_{k_3},\ 3_{(3)}0_{k_3},4_{(3)}0_{k_3},\ 5_{(3)}0_{k_3},\ 6_{(3)}0_{k_3},\ 7_{(3)}0_{k_3},\ 8_{(3)}0_{k_3},\ 9_{(3)}0_{k_3}\},
\\
\, {\cal{N}}_{kd}&=&\{1170_{k_3},\ 1260_{k_3},\ 1350_{k_3},\ 1440_{k_3},\ 2250_{k_3},\ 2340_{k_3},\ 2880_{k_3},\ 4680_{k_3}\},
\\
\, {\cal{N}}_{ke}&=&\{1_{(3)}60_{k_4},\ 11250_{k_4},\ 11340_{k_4},\ 12240_{k_4},\ 12330_{k_4},\ 2_{(3)}30_{k_4},
\nonumber\\
 && 24480_{k_4},\ 22680_{k_4},\ 24660_{k_4},\ 36990_{k_4},\ 4_{(3)}60_{k_4},\ 6_{(3)}90_{k_4}\},
\\
{\cal{N}}_{kf}&=&\{1_{(4)}50_{k_5},\ 1_{(3)}240_{k_5},\ 1_{(3)}330_{k_5},\ 112230_{k_5},\ 12_{(4)}0_{k_5},\ 2_{(3)}480_{k_5},\\
&&\ \ 2_{(3)}660_{k_5},\ 224460_{k_5},\ 24_{(4)}0_{k_5},\ 3_{(3)}990_{k_5},\ 336690_{k_5},\ 36_{(4)}0_{k_5},\ 48_{(4)}0_{k_5}\},
\\
{\cal{N}}_{kg}&=&\{1_{(5)}40_{k_6},\ 1_{(4)}230_{k_6},\ 1_{(3)}2_{(3)}0_{k_6},\ 2_{(5)}80_{k_6},\ 2_{(4)}460_{k_6},\\ 
&&\ \ 2_{(3)}4_{(3)}0_{k_6},\ 3_{(4)}690_{k_6},\ 3_{(3)}6_{(3)}0_{k_6},\ 4_{(3)}8_{(3)}0_{k_6}\},
\end{eqnarray*}
\begin{eqnarray*}
{\cal{N}}_{kh}&=&\{1_{(6)}30_{k_7},\ 1_{(5)}220_{k_7},\ 2_{(6)}60_{k_7},\ 2_{(5)}440_{k_7},\ 3_{(6)}90_{k_7},\ 3_{(5)}660_{k_7},\ 4_{(5)}880_{k_7}
\},
\\
{\cal{N}}_{ki}&=&\{1_{(7)}20_{k_8},\ 2_{(7)}40_{k_8},\ 3_{(7)}60_{k_8},\ 4_{(7)}80_{k_8}\},
\\
{\cal{N}}_{kj}&=&\{1_{(9)}0_{k_9},\ 2_{(9)}0_{k_9},\ 3_{(9)}0_{k_9},\ 4_{(9)}0_{k_9},\ 5_{(9)}0_{k_9},\ 6_{(9)}0_{k_9},\ 7_{(9)}0_{k_9},\ 8_{(9)}0_{k_9},\ 9_{(9)}0_{k_9}\}
\end{eqnarray*}
with the notations $k_i\equiv (k-i),\ i=1,\ \ldots, 9$.
\end{theorem}

\textbf{Proof.}
1. Every element of $\mathcal{S}_k \mathcal{N}_{ka}$ is trivially a PINN because
\[
a0_{(k)} \equiv 0 \pmod{a} \quad \text{for all} \quad a \in \{1, 2, \dots, 9\}.
\]

2. The elements of $\mathcal{S}_k \mathcal{N}_{kb}$ expand as
\begin{eqnarray}
{\cal{S}}_k{\cal{N}}_{kb}&=&\big\{0_{(i)}10_{(j)}20_{(k-2-i-j)},\ 0_{(i)}20_{(j)}10_{(k-2-i-j)},\ 0_{(i)}10_{(j)}80_{(k-2-i-j)},\ 0_{(i)}80_{(j)}10_{(k-2-i-j)},\nonumber\\
&&0_{(i)}20_{(j)}40_{(k-2-i-j)},\ 0_{(i)}40_{(j)}20_{(k-2-i-j)},\ 0_{(i)}20_{(j)}70_{(k-2-i-j)},\ 0_{(i)}70_{(j)}20_{(k-2-i-j)},\nonumber\\
&&0_{(i)}30_{(j)}60_{(k-2-i-j)},\ 0_{(i)}60_{(j)}30_{(k-2-i-j)},\  0_{(i)}40_{(j)}50_{(k-2-i-j)},\ 0_{(i)}50_{(j)}40_{(k-2-i-j)},\nonumber\\
&&0_{(i)}40_{(j)}80_{(k-2-i-j)},\  0_{(i)}80_{(j)}40_{(k-2-i-j)},\ i,j=0,\ 1,\ldots k-2,\ i+j\leq k-2 \big\}.
\end{eqnarray}

We now determine whether each element of $\mathcal{S}_k \mathcal{N}_{kb}$ is a PINN.
$$0_{(i)}10_{(j)}20_{(k-2-i-j)}=10_{(j)}20_{(k-2-i-j)}=3\times(3_{(j)}40_{(k-2-i-j)})=0 \pmod{3},$$
$$0_{(i)}20_{(j)}10_{(k-2-i-j)}=20_{(j)}10_{(k-2-i-j)}=3\times(6_{(j)}70_{(k-2-i-j)})=0 \pmod{3},$$
$$0_{(i)}10_{(j)}80_{(k-2-i-j)}=10_{(j)}80_{(k-2-i-j)}=9\times(1_{(j)}20_{(k-2-i-j)})=0 \pmod{9},$$
$$0_{(i)}80_{(j)}10_{(k-2-i-j)}=80_{(j)}10_{(k-2-i-j)}=9\times(8_{(j)}90_{(k-2-i-j)})=0 \pmod{9},$$
$$0_{(i)}20_{(j)}40_{(k-2-i-j)}=20_{(j)}40_{(k-2-i-j)}=6\times(2_{(j)}40_{(k-2-i-j)})=0 \pmod{6},$$
$$0_{(i)}40_{(j)}20_{(k-2-i-j)}=40_{(j)}20_{(k-2-i-j)}=6\times(6_{(j)}70_{(k-2-i-j)})=0 \pmod{6},$$
$$0_{(i)}20_{(j)}70_{(k-2-i-j)}=20_{(j)}70_{(k-2-i-j)}=9\times(2_{(j)}30_{(k-2-i-j)})=0 \pmod{9},$$
$$0_{(i)}70_{(j)}20_{(k-2-i-j)}=70_{(j)}20_{(k-2-i-j)}=9\times(7_{(j)}80_{(k-2-i-j)})=0 \pmod{9},$$
$$0_{(i)}30_{(j)}60_{(k-2-i-j)}=30_{(j)}60_{(k-2-i-j)}=9\times(3_{(j)}40_{(k-2-i-j)})=0 \pmod{9},$$
$$0_{(i)}60_{(j)}30_{(k-2-i-j)}=60_{(j)}30_{(k-2-i-j)}=9\times(6_{(j)}70_{(k-2-i-j)})=0 \pmod{9},$$
$$0_{(i)}40_{(j)}50_{(k-2-i-j)}=40_{(j)}50_{(k-2-i-j)}=9\times(4_{(j)}50_{(k-2-i-j)})=0 \pmod{9},$$
$$0_{(i)}50_{(j)}40_{(k-2-i-j)}=50_{(j)}40_{(k-2-i-j)}=9\times(5_{(j)}60_{(k-2-i-j)})=0 \pmod{9},$$
$$0_{(i)}40_{(j)}80_{(k-2-i-j)}=40_{(j)}80_{(k-2-i-j)}=12\times(3_{(j)}40_{(k-2-i-j)})=0 \pmod{12},$$
$$0_{(i)}80_{(j)}40_{(k-2-i-j)}=80_{(j)}40_{(k-2-i-j)}=12\times(6_{(j)}70_{(k-2-i-j)})=0 \pmod{12}.$$
Therefore, all elements in $\mathcal{S}_k \mathcal{N}_{kb}$ are PINNs.

3. In the same way, we can prove that all elements of ${\cal{S}}_k{\cal{N}}_{kc}$, ${\cal{S}}_k{\cal{N}}_{kd}$, ${\cal{S}}_k{\cal{N}}_{ke}$, ${\cal{S}}_k{\cal{N}}_{kf}$,  ${\cal{S}}_k{\cal{N}}_{kg}$, ${\cal{S}}_k{\cal{N}}_{kh}$, ${\cal{S}}_k{\cal{N}}_{ki}$ and ${\cal{S}}_k{\cal{N}}_{kj}$ are PINNs. For conciseness, we omit full proof procedures and instead provide a single representative example for each subset.
$${\cal{S}}_k{\cal{N}}_{kc} \ni 0_{(i)}10_{(p)}10_{(q)}10_{(k-3-i-p-q)}=3\times(3_{(p+1)}6_{(q)}70_{(k-3-i-p-q)})=0 \pmod{3},$$
$${\cal{S}}_k{\cal{N}}_{kd} \ni 0_{(i)}10_{(p)}20_{(q)}60_{(k-3-i-p-q)}=9\times(1_{(p+1)}3_{(q)}40_{(k-3-i-p-q)})=0 \pmod{9},$$
$${\cal{S}}_k{\cal{N}}_{ke} \ni 0_{(i)}10_{(p)}20_{(q)}30_{(r)}30_{(k-4-i-p-q-r)}=9\times(1_{(p+1)}3_{(q+1)}6_{(r)}70_{(k-4-i-p-q-r)})=0 \pmod{9},$$
$${\cal{S}}_k{\cal{N}}_{kf} \ni 30_{(p)}30_{(q)}60_{(i)}60_{(k-5-i-p-q)}9=27\times(1_{(p+1)}2_{(q+1)}4_{(i+1)}6_{(k-5-i-p-q)}7)=0 \pmod{27},$$
\begin{eqnarray}
{\cal{S}}_k{\cal{N}}_{kg} &\ni & 10_{(p)}10_{(q)}10_{(i)}10_{(j)}20_{(k-6-i-j-p-q)}3\nonumber\\
&=&9\times(1_{(p+1)}2_{(q+1)}3_{(i+1)}4_{(j+1)}6_{(k-6-i-j-p-q)}7)=0 \pmod{9},\nonumber
\end{eqnarray}
\begin{eqnarray}
{\cal{S}}_k{\cal{N}}_{kh} &\ni & 40_{(p)}40_{(q)}40_{(i)}40_{(j)}40_{(r)}80_{(k-7-i-j-p-q-r)}8\nonumber\\
&=&36\times(1_{(p+1)}2_{(q+1)}3_{(i+1)}4_{(j+1)}5_{(r+1)}7_{(k-7-i-j-p-q)}8)=0 \pmod{36},\nonumber
\end{eqnarray}
\begin{eqnarray}
{\cal{S}}_k{\cal{N}}_{ki} &\ni & 80_{(s)} 40_{(p)}40_{(q)}40_{(i)}40_{(j)}40_{(r)}40_{(k-8-i-j-p-q-r-s)}4\nonumber\\
&=&36\times(2_{(s+1)}3_{(p+1)}4_{(q+1)}5_{(i+1)}6_{(j+1)}7_{(r+1)}8_{(k-8-i-j-p-q-r-s)}9)=0 \pmod{36},\nonumber
\end{eqnarray}
\begin{eqnarray}
{\cal{S}}_k{\cal{N}}_{kj} &\ni & 10_{(s)} 10_{(p)}10_{(q)}10_{(i)}10_{(j)}10_{(r)}10_{(m)}10_{(k-9-i-j-p-q-r-s-m)}1\nonumber\\
&=&9\times(1_{(s+1)}2_{(p+1)}3_{(q+1)}4_{(i+1)}5_{(j+1)}6_{(r+1)}7_{(m+1)}8_{(k-9-i-j-p-q-r-s)}9)=0 \pmod{9}.\nonumber
\end{eqnarray}
The theorem is proved. $\square$

From Theorem~\ref{mainthm} and convention~\eqref{convention}, we conclude that
\begin{equation}
    \mathrm{PINN}_{k+1} \supset \mathrm{PINN}_k.
\label{kck1}
\end{equation}

\section{Further Properties and Conditions}

\subsection{Zero Insertion in Repdigit PINNs}

We now address the question: Can zeros be inserted into other repeated--digit PINNs from \eqref{3na} while preserving the Niven property? Regrettably, the answer appears negative. Consider the PINN $1_{(27)} \equiv 0 \pmod{27}$. If we insert a zero to form $1_{(26)}01$, the resulting number is not an NN (and consequently not a PINN) since
\[
1_{(26)}01 \equiv 18 \not\equiv 0 \pmod{27}.
\]
Similarly, zero insertion fails for other PINNs in \eqref{3na}. For instance:
\begin{itemize}
  \item For PINN $1_{(81)} \equiv 0 \pmod{81}$, $1_{(80)}01 \equiv 72 \not\equiv 0 \pmod{81}$.
  \item For PINN $1_{(111)} \equiv 0 \pmod{111}$, $1_{(110)}01 \equiv 102 \not\equiv 0 \pmod{111}$.
  \item With double zero insertion, $1_{(110)}001 \equiv 12 \not\equiv 0 \pmod{111}$.
\end{itemize}
Thus, $1_{(80)}01$, $1_{(110)}01$, and $1_{(110)}001$ are neither NNs nor PINNs.

\subsection{Symmetry Property}

All $k$--digit PINNs exhibit $\mathcal{S}_k$ symmetry:
\begin{equation}
\mathcal{S}_k \mathrm{PINN}_k = \mathrm{PINN}_k.
\label{invk}
\end{equation}

\subsection{Necessary and Sufficient Conditions}

For the seventh problem, the necessary and sufficient conditions for a $k$--digit PINN, denoted $A_k = a_k a_{k-1} \ldots a_2 a_1$, follow directly from its definition:
\begin{equation}
s_j \sum_{i=1}^k a_i \times 10^{i-1} = b_j \sum_{i=1}^k a_i, \quad j = 1, 2, \ldots, k!,\  s_j \in \mathcal{S}_j,
\label{sn}
\end{equation}
where $a_i \in \{0, 1, 2, \ldots, 9\}$ for $i = 1, 2, \ldots, k$, and $b_j$ (for $j = 1, 2, \ldots, k!$) are arbitrary possible positive integers.

Although the PINNs given in Theorem~\ref{mainthm} \eqref{PINNk} satisfy \eqref{sn}, solving this system becomes prohibitively difficult for large $k$. The two-step search procedure employed in this paper provides a more efficient method for finding PINNs.

\subsection{Digit Sum Conditions}

For the eighth problem, Theorem~\ref{mainthm} implies that trivial PINNs with a single nonzero digit of the form $a0_{(0)}$ (where $a \in \mathrm{NN}_1$) have digit sum $a$. All other PINNs have digit sums divisible by 3. Excluding the repdigit PINNs in \eqref{3na} and the trivial cases $a0_{(0)}$ ($a \in \mathrm{NN}_1$), we have only found PINNs satisfying $3 \leq \sum(\mathrm{PINN}) \leq 81$. We conjecture that no exceptions exist beyond these cases.

\subsection{Relations to Other Sequences}

For the ninth problem, it is evident that PINNs constitute a subset of NNs. Repdigit NNs are naturally PINNs. We have not yet identified any relationships between PINNs and established OEIS sequences. Apart from repdigit PINNs, we have not discovered any other distinguished subclasses such as palindromic numbers or prime numbers.

From the repdigit PINNs presented in \eqref{3na}, we can extract an interesting subset of primes:
\begin{eqnarray}
&&\{3,\, 37,\, 163,\, 757,\, 1999,\ 8803,\ 9397,\, 13627,\, 15649,\, 231643,\, 313471,\, 333667,\, 338293,\ 1014877,\nonumber\\
&&\quad 1056241,\  1168711,\, 2028119,\, 2064529,\, 2462401,\, 2558791,\, 4448359,\, 9438277,\, 34720813,\, \nonumber\\
&&\quad 86455449,\, 104620573,\, 127020961,\, 178064569,\, 247629013,\, 618846643,\ 440334654777631,\nonumber\\
&&\quad 676421558270641,\,\, 2212394296770203368013,\,\, 130654897808007778425046117,\, \ldots\ \}.
\end{eqnarray}
Some of these primes appear in \eqref{knmpqr}, while others may also belong to known prime collections. Further investigation of primes related to repdigit PINNs extends beyond this paper's scope and will be addressed in subsequent work.

\section{Open Problems and Future Research}

For the final discussion, we present several significant open questions for future research:
\begin{enumerate}
    \item Can the concept of PINNs be extended to arbitrary number bases?

    \item How can we exhaustively characterize all repdigit PINNs?

    \item Is there a systematic approach for solving algebraic equations of the form \eqref{0} using distinguished prime sets?

    \item Are there more efficient algorithms for generating high--digit PINNs?

    \item To what extent do permutation invariant properties hold for other variants of extended Niven numbers,
such as $c$--NNs, higher--order NNs, multiple NNs, and generalized NNs?
\end{enumerate}

\section*{Acknowledgement} {\noindent \small
 The work was sponsored by the National Natural Science Foundations of China (Nos.12235007, 12375003, 11975131).
 The authors are indebted to thank Profs. M. Jia, Q. P. Liu, X. B. Hu, B. F. Feng and R. X. Yao for their helpful discussions.}


\begin{thebibliography}{99}
\bibitem{Niven}Niven I. 1977, Integers divisible by the sum of their digits, In: Mathematics of Computation, 1943-1993: A Half--Century of Computational Mathematics (W. Gautschi, ed.), Proceedings of Symposia in Applied Mathematics, Vol. 48. American Mathematical Society.
\bibitem{Kennedy82}Kennedy R. E. 1982, Digital sums, Niven numbers, and natural density, Fibonacci Quarterly, 20, 160--162.
\bibitem{Cooper} Cooper C. N. and Kennedy R. E. 1993, On Consecutive Niven Numbers, Fibonacci Quarterly, 31, 146--151.
\bibitem{Grundman}Grundman H. G. 1994, Sequences of Consecutive Niven Numbers, Fibonacci Quarterly, 32, 174--175.
\bibitem{Koninck}De Koninck J. M. and Doyon N. 2003, On the Number of Niven Numbers up to \(x\),
   Fibonacci Quarterly, 41, 431--440.
\bibitem{Cai}Cai T. 2003, On 2--Niven Numbers and 3--Niven Numbers, Fibonacci Quarterly, 41, 168--173.
\bibitem{Fredricksen}Fredricksen H., Ionascu E. J. and Luca F. 2008, Minimal Niven Numbers,
   INTEGERS: The Elect. J. Combin. Number Theory, 8, A19.
\bibitem{Sanna}Sanna C. 2021, Additive bases and Niven numbers, Bulletin  Australian Math. Soc., 104, 373--380
\bibitem{Wilson}Wilson B. 1996, Construction of small consecutive Niven numbers, Fibonacci Quarterly,
34, 240--43.
\bibitem{Olver}Olver P. J., 1993, Applications of Lie groups to differential equations, Springer, New York.
\bibitem{Lou}Lou S. Y., 2022, Nonlocal symmetries of nonlinearintegrable systems, In: Nonlinear systems and their remarkable mathematical sturctures, Vol. 3: Contribiutions from China, Eds: N. Euler and D. J. Zhang, First edition published by CRC Press 6000 Broken Sound Parkway NW, Suite 300, Boca Raton, FL 33487--2742 and by CRC Press2 Park Square, Milton Park, Abingdon, Oxon, OX144RN.
\bibitem{Lou1}Lou S. Y. and Feng B. F., 2025, Symmeties and integrable systems, Fundamental Research, 5, 1947-1953.
\end{thebibliography}
\end{document}